\documentclass[10pt]{amsart}
\usepackage{graphicx}
\usepackage[english,francais]{babel}
\usepackage{color}
\usepackage[ansinew]{inputenc}
\usepackage[active]{srcltx}
\usepackage[all]{xy}
\usepackage{amssymb}
\usepackage{amsmath}

\newtheorem{thm}{Th\'{e}or\`{e}me}[section]

\theoremstyle{definition}

\newtheorem{definition}[thm]{D\'{e}finition}

\theoremstyle{remark}

\numberwithin{equation}{section}

\newcommand{\ra}{\rightarrow}

\newcommand{\Pun}{\mathbb P^{1}_{\mathbb C}}

\newcommand{\Gm}{\mathbb G_{m}}
\newcommand{\mgg}{\mathcal M_{\mathbb G_{m}}^{\mathbb G_{m}}}

\title{Fibre de Milnor motivique \`a l'infini }

\author{ Michel Raibaut }

\address{Laboratoire JA Dieudonn\'e, UMR CNRS 6621\\ Universit\'e Nice Sophia-Antipolis, parc Valrose, 06108 Nice Cedex 2, France}

\email{raibaut@unice.fr}

\begin{document}
\maketitle

Note re\c{c}ue par le Compte Rendu de l'Acad\'emie des Sciences le 1 D\'ecembre 2009 et accept\'ee apr\`es r\'evision le 13 Janvier 2010.\\

\begin{abstract}
\selectlanguage{francais}
Pour une application r\'eguli\`ere $f:U \ra \mathbb A^{1}$ \`a source lisse, nous d\'efinissons une fibre de Milnor motivique \`a l'infini et nous la calculons dans le cas d'un polyn\^ome de Laurent non d\'eg\'en\'er\'e pour son poly\`edre de Newton \`a l'infini.

\vskip 0.5\baselineskip

\selectlanguage{english}
\noindent{ Abstract.}
Given a regular map $f:U \ra \mathbb A^{1}$ on a smooth variety $U$ we define a motivic Milnor fiber at infinity and we compute it in the case of a non degenerate Laurent polynomial for its Newton polyhedra at infinity.
\end{abstract}

\selectlanguage{francais}

\section{Introduction}
Soit $U$ une vari\'et\'e complexe lisse et $f:U\ra \mathbb A^{1}$ une application r\'eguli\`ere non constante.
 Il existe $R>0$ tel que $f: U \setminus f^{-1}(\overline{D}(0,R)) \ra \mathbb C
\setminus \overline{D}(0,R)$ est une fibration topologique localement triviale \cite{Pham}.
Les espaces de cohomologie \`a support compact
$H_{c}^{*}\left(f^{-1}(t),\mathbb Q\right)$ de la fibre en $t$ sont munis d'une structure de Hodge mixte. Steenbrink et Zucker  puis M. Saito ont montr\'e comment construire une structure de Hodge mixte limite lorsque $t$ tend vers l'infini. Sabbah \cite{Sabbah} l'a retrouv\'ee en consid\'erant la transformation de Fourier sur des modules convenables sur l'anneau des op\'erateurs diff\'erentiels. Le spectre de cette structure limite  est un invariant de $f$ appel\'e \textbf{spectre \`a l'infini} \cite{Sabbah} (5.4), \cite{GLN} et \cite{NS}.\\
Pour un morphisme $f:X\ra \mathbb A^{1}_{\mathbb C}$ avec $X$ lisse et $x\in f^{-1}(0)$, Denef et Loeser \cite{DLBarcelona} et \cite{DLIguza} obtiennent le spectre de Hodge-Steenbrink de $f$ en $x$ \`a partir de la fibre de Milnor motivique $S_{f,x}$. Guibert, Loeser et Merle \cite{GLM1} g\'en\'eralisent $S_{f,x}$ en construisant une fibre de Milnor motivique $S_{f,U}$ adapt\'ee \`a un ouvert $U$ de $X$.
Ainsi, en regardant la vari\'et\'e $U$ comme ouvert dans une compactification, nous reconsid\'erons le probl\`eme initial du point de vue motivique et d\'efinissons  une \textbf{fibre de Milnor motivique \`a l'infini }
$\boldsymbol{S_{f,\infty}}$ . C'est un invariant appartenant \`a $\mgg$ (anneau de Grothendieck des vari\'et\'es au dessus de $\Gm$ munies d'une action de $\Gm$ \cite{GLM1}, \cite{GLM2}). Il ne d\'epend pas de la compactification choisie (\ref{independance}). Dans le groupe de Grothendieck $K_{0}(SH^{mon})$ des structures de Hodge munies d'un automorphisme d'ordre fini, il se r\'ealise en la classe de la structure de Hodge mixte \`a l'infini de $f$ (\ref{compatibilite}). Il donne ainsi acc\`es au spectre \`a l'infini de $f$ (\ref{thmspectre}). Nous le calculons pour un polyn\^ome de Laurent non d\'eg\'en\'er\'e pour son poly\`edre de Newton \`a l'infini (\ref{calcul}).
\section{Fibre de Milnor motivique \`a l'infini} 
Pour les d\'efinitions usuelles des anneaux de Grothendieck, des espaces d'arcs et de la mesure motivique $\mu$, on pourra se r\'ef\'erer \`a \cite{DLBarcelona}, \cite{GLM1}
et \cite{GLM2}. On appelle vari\'et\'e tout $\mathbb C$-sch\'ema s\'epar\'e r\'eduit de type fini sur $\mathbb C$.
\noindent Nous consid\'erons dans la suite une vari\'et\'e lisse $U$, une fonction r\'eguli\`ere $f : U \ra \mathbb A^{1}$, l'ouvert $U^{*}=U\setminus f^{-1}(0)$ et l'immersion $j:\mathbb A^{1}_{\mathbb C} \ra \mathbb P^{1}_{\mathbb C},\:a \mapsto [1:a]$. \\ 

\begin{definition}
On appelle \textbf{compactification de $f$} tout triplet $(X,\:i_{X}:  U \ra X,\:f_{X}:X \ra \mathbb P^{1})$ o\`u X est une vari\'et\'e, $i_{X}$ est une immersion ouverte dominante, $f_{X}$ est une application propre et $ f_{X} \circ i_{X}=j \circ f$.\\
\end{definition}
Soit $(X,i,\hat{f})$ une compactification de $f$, notons $F$ le ferm\'e $X\setminus i(U^{*})$, $\hat{f}$ la fonction 
$ [{\hat{f}}_{(0)}:{\hat{f}}_{(1)}]$ et $\hat{f}^{(\infty)}$ la fonction $\frac{{\hat{f}}_{(0)}}{{\hat{f}}_{(1)}} : X\setminus X_{0} \ra \mathbb A^{1}_{\mathbb C}$ o\`u  $X_{0}$ et $X_{\infty}$ sont les vari\'et\'es $\hat{f}_{(1)}^{-1}(0)$ et $\hat{f}_{(0)}^{-1}(0)$.
Nous travaillons dans l'anneau de Grothendieck $\mathcal M^{\Gm}_{X_{\infty}\times \Gm}$ ((2.2) \cite{GLM1} et (2.2) \cite{GLM2}) des vari\'et\'es $\left(V\overset{(p_{X},p_{\Gm})}{\longrightarrow} X_{\infty} \times \Gm,\sigma \right)$ o\`u $\sigma$ est une bonne action sur $V$ du groupe multiplicatif $\Gm$, $p_{X}$ est un morphisme \`a fibres $\Gm$-invariantes et $p_{\Gm}$ est un morphisme homog\`ene pour l'action $\sigma$.\\ 
Pour  $(n, \delta)\in \mathbb N^{*2}$ on pose  $X_{n}^{\delta}:=\left\{ \varphi \in \mathcal L(X\setminus X_{0}) \mid ord_{t}\hat{f}^{(\infty)}(\varphi)=n, ord_{t}\varphi^{*}(\mathcal I_{F})\leq \delta n \right\}$. C'est une partie semi-alg\'ebrique de l'espace des arcs $\mathcal L(X)$ de $X$ munie de l'action usuelle de $\Gm$ sur les arcs $\lambda.\varphi(t)=\varphi(\lambda t)$ et du morphisme $\varphi \mapsto \left(\varphi(0),ac\left(\hat{f}^{(\infty)}(\varphi)\right)\right)$ vers $X_{\infty}\times \Gm$, o\`u  $ac\left(\hat{f}^{(\infty)}(\varphi)\right)$ est le premier coefficient non nul de la s\'erie $\hat{f}^{(\infty)}(\varphi(t))$.
Consid\'erons alors la fonction z\^eta motivique $Z_{\hat{f}^{(\infty)},i(U^{*})}^{\delta}(T)$ ((3.7) \cite{GLM1}), c'est la s\'erie g\'en\'eratrice de la suite des mesures  motiviques 
$\mu(X_{n}^{\delta})$ d\'efinie par 
$$Z_{\hat{f}^{(\infty)},i(U^{*})}^{\delta}(T):=\underset{n \geq 1}{\sum}\:\mu(X_{n}^{\delta})T^{n}\in \mathcal M^{\Gm}_{X_{\infty}\times \Gm}[[T]].$$
Si la compactification est lisse alors cette fonction z\^eta est rationnelle et $S_{\hat{f}^{(\infty)},i(U^{*})}:=-\underset{T \ra \infty}{lim}Z_{\hat{f}^{(\infty)},i(U^{*})}^{\delta}(T)$ \'el\'ement de  $\mathcal M_{X_{\infty} \times \Gm}^{\Gm}$ ne d\'epend pas de $\delta$ pour $\delta$ assez grand ((3.8) \cite{GLM1}). Si la compactification est singuli\`ere, le lieu singulier est alors contenu dans le ferm\'e $X\setminus i(U)$. La m\^eme preuve fonctionne et le r\'esultat est inchang\'e.
Pour une vari\'et\'e $S$, notons
$p_{S!}:\mathcal M_{S \times \Gm}^{\Gm} \ra \mathcal M_{\Gm}^{\Gm}, (V \overset{p,q}{\ra} S \times \Gm, \sigma) \mapsto 
(V \overset{q}{\ra} \Gm, \sigma) $. Ainsi, \\
\begin{thm}\label{independance} 
Soit $(X,i_{X},f_{X})$ et $(Y,i_{Y},f_{Y})$ deux compactifications de $f$. On a 
$$p_{X_{\infty}!}S_{f_{X}^{(\infty)},i_{X}(U^{*})}=p_{Y_{\infty}!}S_{f_{Y}^{(\infty)},i_{Y}(U^{*})}\in \mathcal M_{\mathbb G_{m}}^{\mathbb G_{m}}.$$
On note cette valeur $S_{f, \infty}\in \mathcal M_{\mathbb G_{m}}^{\mathbb G_{m}} $ et on l'appelle \textbf{fibre de Milnor motivique \`a l'infini de $f$} au sens de (\cite{DLBarcelona}, \cite{GLM1}).
C'est un nouvel invariant de $f$ que l'on peut tenter de calculer \`a partir de n'importe quelle compactification lisse ou non de $f$.\\ \\
\end{thm}

\section{Lien avec la structure de Hodge mixte limite.}
Soit $(X,i,\hat{f})$ une compactification de $f$. Par ((3.9)\cite{GLM1}), il existe un unique morphisme de $\mathcal M_{\mathbb C}$-modules  $S_{\hat{f}^{\infty}}:\mathcal M_{X\setminus X_{0}} \ra \mathcal M_{X_{\infty}\times \Gm}^{\Gm}$ tel que pour tout 
morphisme propre $p:Z \ra X\setminus X_{0}$ avec $Z$ lisse et pour tout ouvert dense $V$ de $Z$,
$S_{\hat{f}^{\infty}}([V\ra X\setminus X_{0}])$ vaut $p_{!}(S_{\hat{f}^{\infty}\circ p,V})$. On note $MHM_{X\setminus X_{0}}$ la cat\'egorie ab\'elienne des modules de Hodge mixtes sur $X\setminus X_{0}$, $K_{0}(MHM_{X\setminus X_{0}})$ l'anneau de Grothendieck correspondant vu comme $\mathcal M_{\mathbb C}$-module et $\Psi_{\hat{f}^{\infty}}$ le foncteur cycles proches \cite{S}. Par additivit\'e il existe un unique morphisme $\mathcal M_{\mathbb C}$-lin\'eaire $H:\mathcal M_{X\setminus X_{0}} \ra K_{0}(MHM_{X\setminus X_{0}})$ 
tel que pour tout $p:Z \ra X\setminus X_{0}$ avec $Z$ lisse, $H([p:Z\ra X\setminus X_{0}])$ est la classe $[Rp_{!}(\mathbf Q_{Z})] $ o\`u $\mathbf Q_{Z}$ est le module de Hodge trivial sur $Z$ (lemme (14.61) \cite{PS}). On construit de m\^eme un morphisme $H:\mathcal M_{X_{\infty}\times \Gm}^{\Gm} \ra K_{0}(MHM_{X_{\infty}}^{mon})$ ((3.16)\cite{GLM1}) compatible avec l'action de $\Gm$ et la monodromie.
En appliquant ((3.17) \cite{GLM1}) et la compatibilit\'e des modules de Hodge mixtes avec l'image directe 
$p: X_{\infty} \ra Spec(\mathbb C)$ on a avec les notations ci dessus : \\

\begin{thm} \label{compatibilite} Pour une compactification $(X,i,\hat{f})$ le diagramme suivant est commutatif :
\tiny{
$$\xymatrix{
    \mathcal M_{X\setminus X_{0}} \ar[r]^-{S_{\hat{f}^{\infty} }} \ar[d]_-{H} 
   & \mathcal M_{X_{\infty} \times \Gm}^{\Gm}  \ar[r]^-{ p_{!} } \ar[d]_-{H} 
   & \mgg  \ar[d]_-{H} 
     \\
   K_{0}\left(MHM_{X \setminus X_{0}}\right) \ar[r]^-{\Psi_{\hat{f}^{\infty}} }        
   & K_{0}\left(MHM_{X_{\infty}}^{mon}\right) \ar[r]^-{p_{!}}        
   & K_{0}\left(MHM_{Spec \mathbb C}^{mon}\right) 
     \\
  }$$
}\normalsize{}
En particulier $H(S_{f,\infty})=p_{!}(\Psi_{\hat f^{\infty}}(Ri_{!}\mathbf Q_{U}))$.\\
\end{thm}
La structure de Hodge mixte limite sur  $H^{k}_{c}\left( f^{-1}(t),\mathbb Q\right)$ s'identifie \cite{Sabbah} \`a la structure de Hodge mixte de 
$\mathbb H_{c}^{k} (\hat{f}^{-1}(\infty), \psi_{1/f}(Ri_{!}\underline{\mathbb Q}_{U}))$. On obtient ce groupe d'hypercohomologie en prenant le faisceau pervers sous jacent de
$p_{!}(\Psi_{\hat f^{\infty}}(Ri_{!}\mathbf Q_{U})) \in D^{b}\left(MHM_{Spec \mathbb C}^{mon}\right)$   ((14.1.1)\cite{PS} et \cite{S}). Notons $\Phi$ le morphisme $K_{0}\left(MHM_{Spec \mathbb C}^{mon}\right) \ra K_{0}(SH^{mon})$ ((6.1)\cite{GLM1} et (14.1.1)\cite{PS}). Par d\'efinition  ((3.1.2) \cite{DLBarcelona}) la classe de la structure de Hodge mixte limite est $\underset{k}{\sum} (-1)^{k}[ \mathbb H_{c}^{k} (\hat{f}^{-1}(\infty), \psi_{\frac{1}{f}}(Ri_{!}\underline{\mathbb Q}_{U}))]$. Elle est \'egale \`a $\Phi(p_{!}(\Psi_{\hat f^{\infty}}(Ri_{!}\mathbf Q_{U})))$, elle m\^eme \'egale \`a $\Phi(H(S_{f,\infty}))$ par (\ref{compatibilite}). Le spectre \`a l'infini est le spectre de la classe de la structure de Hodge mixte limite donn\'e par le spectre de Hodge $sph : K_{0}\left(SH^{mon}\right) \ra \mathbb Z[\mathbb Q]$ ((6.1.4) \cite{GLM1}). Ainsi, \\
\begin{thm} \label{thmspectre} 
Pour un morphisme non constant $f : U \ra \mathbb A^{1}$ avec $U$ lisse, la classe de la structure de Hodge mixte limite \`a l'infini est $\Phi(H(S_{f,\infty}))$ et 
son spectre \`a l'infini vaut $Sp(S_{f,\infty})$ o\`u $Sp = sph \circ \Phi \circ H$.
\end{thm}

\section{Calcul dans le cas d'un polyn\^ome de Laurent non d\'eg\'en\'er\'e.}
Soit  $f(\underline{x}):=\underset{\alpha \in \mathbb Z^{d}}{\sum}a_{\alpha}\underline{x}^{\alpha} \in \mathbb C[x_{1},..,x_{d}][x_{1}^{-1},..,x_{d}^{-1}]$, son support $\text{supp}(f)$ est l'ensemble $\{\alpha \in \mathbb N^{d} \mid a_{\alpha}\neq 0\}$ et son poly\`edre de Newton \`a l'infini $\Gamma_{-}$ est l'enveloppe convexe de $\text{supp}(f)\cup\{0\}$. On note $\Gamma$ les faces de $\Gamma_{-}$ ne contenant pas l'origine.
Pour toute face $\gamma$ de $\Gamma_{-}$, on note $f_{\gamma}(\underline{x})$ le polyn\^ome quasi-homog\`ene $\underset{\alpha \in \gamma}{\sum}a_{\alpha}\underline{x}^{\alpha}$. Au sens de Kouchnirenko \cite{K}, supposons que $f$ est  \textbf{non d\'eg\'en\'er\'e pour son poly\`edre de Newton \`a l'infini} : pour toute face $\gamma$ de $\Gamma$, le polyn\^ome  $f_{\gamma}$ est lisse sur $\Gm^{d}$. 
Dans $ \mathbb P^{d}\times \mathbb P^{1}$, on choisit la compactification 
$X= \{([x],[\alpha:\beta]) \mid \alpha\overset{\sim}{P}(x)x_{0}^{deg(Q)}=\beta x_{0}^{deg(P)}Q(x) \},$  avec $\hat{f} :  X  \ra  \Pun$, $([\underline{x}],[\alpha:\beta]) \mapsto [\alpha : \beta]$  et $i : \mathbb A_{\mathbb C}^{d}  \ra  X$, $x  \mapsto  ([1:x], [f(x),1])$ o\`u $f=\frac{P}{Q}$ avec $Q$ monomial, $P$ et $Q$ premiers entre eux et $\tilde{P}$ l'homog\'en\'eis\'e de $P$.
Comme dans \cite{Gil}, on prouve la rationalit\'e de 
$Z^{\delta}_{\hat{f}^{\infty},i(U^{*})}(T)$ sans recours au th\'eor\`eme de r\'esolution des singularit\'es d'Hironaka.
Dans l'espace affine $U^{*}=\mathbb G_{m}^{d}\setminus f^{-1}(0)$, on utilise des arcs de Laurent $\varphi=\left( \frac{P_{i}(t)}{t^{\omega_{i}}} \right)$ dont l'origine appartient \`a $X_{\infty}$ et 
v\'erifiant la condition au bord $ord_{t}  \varphi^{*}(\mathcal{I}_{F}) \leq \delta\; ord_{t} \hat{f}^{\infty} \varphi$. 
Pour ces arcs, $P_{i}(t)$ est une s\'erie formelle inversible et $(\omega_{i})$ appartient \`a 
$\Omega  =\{\omega \in \mathbb Z^{d} \mid  max((\omega\mid.)_{\mid\Gamma_{-}})>0\}$. Pour chaque arc, on note $\gamma(\omega)$ la face de $\Gamma_{-}$ o\`u 
la forme lin\'eaire $(\omega\mid.)_{\mid\Gamma_{-}}$ atteint son maximum. Cette face ne contient pas $0$. Les arcs sont donc class\'es par les faces de $\Gamma$ et pour toute face $\gamma \in \Gamma$ on note
$C_{\gamma}:=\left\{\omega \in \Omega \mid \gamma(\omega)=\gamma \right\} $.
La non d\'eg\'enerescence de $f$ correspond \`a la lissit\'e des $f_{\gamma}$ et permet de mesurer les espaces d'arcs. On obtient alors, \\
\begin{thm} \label{calcul}
Soit  $f\in \mathbb C[x_{1},..,x_{d}][x_{1}^{-1},..,x_{d}^{-1}]$
non d\'eg\'en\'er\'e pour son poly\`edre de Newton $\Gamma$. La fibre de Milnor motivique \`a l'infini de $f : \Gm^{d} \ra \mathbb A^{1}$ vaut \\
$$S_{f,\infty}= -\underset{\gamma \in \Gamma}{\sum}\:\:
\chi(C_{\gamma}) \left[\Gm^{d}\setminus f_{\gamma}^{-1}(0),\:f_{\gamma}^{-1},\:\sigma(\gamma ) \right] $$
o\`u $\chi$ est la caract\'eristique d'Euler \`a support compact  
et $\sigma\left(\gamma \right)$ est une action de $\Gm$ sur $\Gm^{d}\setminus f_{\gamma}^{-1}(0) $  de la forme
$\sigma\left(\gamma \right)(\lambda,x)=(\lambda^{-\omega_{i}}x_{i})$ avec $\omega \in C_{\gamma}$. La classe 
$\left[\Gm^{d}\setminus f_{\gamma}^{-1}(0),\:f_{\gamma}^{-1},\:\sigma(\gamma ) \right]$ ne d\'epend pas de $\omega$.\\ \\
Si de plus le polyn\^ome est commode (0 est contenu dans l'int\'erieur du poly\`edre de Newton $\Gamma_{-}$) alors $\chi(C_{\gamma})$ est nulle
pour toute face  $\gamma$ contenue dans un hyperplan de coordonn\'ees  et vaut $(-1)^{d-dim(\gamma)}$ sinon.\\ \\
Par (\ref{thmspectre}), le spectre \`a l'infini de $f$ vaut: 
$sp(f)=-\underset{\gamma\in \Gamma}{\sum}\: \chi(C_{\gamma}) \:
Sp\left [f_{\gamma}^{-1}(1),\mu_{\gamma}\right]$ \\
o\`u $\left(f_{\gamma}^{-1}(1),\mu_{\gamma}\right)$ est la vari\'et\'e $f_{\gamma}^{-1}(1)$ munie de l'action de  induite par $\sigma\left(\gamma\right)$.\\
\end{thm}

Pour $f\in \mathbb C[x_{1},..,x_{d}]$ non d\'eg\'en\'er\'e pour son poly\`edre de Newton on obtient une formule similaire. Pour cela, comme dans \cite{GLM2}, on stratifie $\mathbb A^{d}$ en produit de tores et on utilise l'additivit\'e de la fibre de Milnor motivique ((3.9)\cite{GLM1}). On applique alors le th\'eor\`eme (\ref{calcul}) \`a la restriction de $f$ \`a chaque strate.

\end{document}